\documentclass[final]{siamltex}

\divide \time 6 \newcommand{\timestamp}{\the\year-\the\month-\the\day-\the\time}

\usepackage{amsmath,amssymb,amsfonts,latexsym, mathrsfs}

\newlength{\tailwidth}
\newlength{\xxxwidth}

\usepackage[dvipsone]{epsfig}

\usepackage{color}

\newtheorem{algorithm}[theorem]{Algorithm}

\usepackage{alltt}
\newcommand{\mc}{\multicolumn}

\newcommand{\mcc}[1]{\mc{1}{|c|}{#1}}
\newcommand{\importfig}[2][3.0in]{%
   \begin{minipage}{#1}
   \begin{tabular}{@{}c@{}}
   \ifthenelse{\equal{#1}{!}}%
   {\psfig{file=#2}}%
   {\psfig{
      file=#2,width=#1}
      \\%
    {\begin{minipage}{#1}\begin{alltt}#2\end{alltt}\end{minipage}}%
   }%
   \end{tabular}
   \end{minipage}
   }

\title{ML($n$)BiCGStabt: A ML($n$)BiCGStab Variant with ${\bf A}$-transpose
}

\author{
Man-Chung Yeung\thanks{Dept. 3036, 1000 East University Avenue,
Laramie, WY 82071. E-mail: myeung@uwyo.edu.
}
}

\begin{document}
\setcounter{page}{1}

\maketitle

\begin{abstract}
The 1980 IDR method\cite{ws} plays an important role in the history of Krylov subspace methods. It started the research of transpose-free Krylov subspace methods. In this paper, we make a first
attempt to bring back ${\bf A}$-transpose to the research area by presenting a new ML($n$)BiCGStab variant
that involves
${\bf A}$-transpose in its implementation.
Comparisons of this new algorithm with the existing ML($n$)BiCGStab algorithms will be presented.

\end{abstract}



\pagestyle{myheadings} \thispagestyle{plain} \markboth{M. YEUNG
}{ML($n$)BiCGStabt}

\section{Introduction}
\label{Introduction}
ML($n$)BiCGStab is a transpose-free Krylov subspace method for the solution of linear systems
\begin{equation} \label{equ:7-6-1}
{\bf A} {\bf x} = {\bf b} 
\end{equation}
where ${\bf A} \in {\mathbb
 C}^{N \times N}$
and ${\bf b} \in {\mathbb C}^N$. It was introduced by Yeung and Chan\cite{yeungchan} in 1999 and its algorithms were recently reformulated
by Yeung\cite{ye2012}. ML($n$)BiCGStab is a natural generalization of BiCGStab\cite{van},
built from a multiple starting BiCG-like algorithm called ML($n$)BiCG, through the Sonneveld-van der Vorst-Lanczos procedure (SVLP), namely,
the procedure introduced by Sonneveld\cite{sonn} and van der Vorst\cite{van} in the construction of CGS and BiCGStab from BiCG\cite{fletcher}.
In theory, ML($n$)BiCGStab is a method that lies between the Lanczos-based BiCGStab and the Arnoldi-based GMRES/FOM\cite{saad2}. In fact,
it
is a BiCGStab when $n = 1$ and becomes a GMRES/FOM when $n = N$ (see \cite{ye2012, yeung7}).
In computation, ML($n$)BiCGStab can be much more stable and converge much faster than BiCGStab.
We once tested
it
on the standard oil reservoir simulation test data called SPE9 which contains a sequence of linear systems and found that it
reduced
the total computational time by 
$60\%$ when compared to BiCGStab. Tests made on the data from matrix markets
also supported the superiority of ML($n$)BiCGStab over BiCGStab. For details, one is referred to \cite{ye2012, yeungchan}.

The author once constructed a new version of ML(n)BiCG where the left residuals are not just given by the monomial basis,
but are orthogonalized 
against previous right-hand side
residuals. In structure, this new ML($n$)BiCG is closer to the classical BiCG than the one in \cite{yeungchan} is.
Numerical experiments, however, showed that this new ML(n)BiCG was
unstable and
weaker than the standard BiCG.
Moreover, in \cite{yeungboley}, Yeung and Boley derived a SVLP from a one-sided multiple starting band Lanczos procedure (MSLP) with $n$ left-starting and $m$ right-starting vectors respectively. From their experiments with multi-input multi-output time-invariant linear dynamical systems, they observed that SVLP is more stable than MSLP when $m \ne n$.
The
two examples of comparison
hint that, when $m \ne n$,
a stable
multiple starting
procedure
with $\bf A$-transpose
may come
from
a modification of a SVLP. In this paper, we make a first step in this direction by introducing $\bf A$-transpose into ML($n$)BiCGStab. We call the resulting algorithm ML($n$)BiCGStabt, standing for ML($n$)BiCGStab with transpose.



There exist two ML($n$)BiCGStab algorithms, labeled as Algorithms 4.1 and 5.1 respectively in \cite{ye2012}, derived from different definitions of the residual vectors ${\bf r}_k$.
While both algorithms are numerically stable in general,
one is relatively more stable than the other.
ML($n$)BiCGStabt is a modified version of Algorithm 5.1 so that it enjoys the same
level of
stability with Algorithm 4.1.

Other
extensions
of IDR, CGS and BiCGStab exist. Among them are BiCGStab2\cite{gut}, 
BiCGStab($l$)\cite{SF}, 
GPBi-CG\cite{zhang}, 
IDR($s$)\cite{
sonn1, gs2013}, 
IDRstab\cite{sv2013}, and GBi-CGSTAB($s,l$)\cite{tasu}. Related articles include \cite{dsyyz, gut1, gutzem, sonn2013}.

The outline of the paper is as follows.
In \S\ref{sec:8-1}, index functions in \cite{yeungboley} are introduced. They are helpful in the
construction of a ML($n$)BiCGStab algorithm.
In \S\ref{sec:12-24-1}, we present the
ML($n$)BiCG algorithm from \cite{yeungchan}. The derivation of every ML($n$)BiCGStab algorithm is based on it.
In \S\ref{sec:12-24-2}, we introduce the ML($n$)BiCGStabt algorithm
and its properties.
In \S\ref{sec:12-19}, numerical experiments are presented,
and in \S\ref{con}, concluding remarks are given.

\section{Index Functions}\label{sec:8-1} Let be given a $n \in {\mathbb N}$, the set of positive integers.
For all $k \in {\mathbb Z}$, the set of all integers, we define
$$
\begin{array}{lcl}
g_n(k) = \lfloor ( k - 1 ) /n \rfloor & \mbox{and}& r_n(k) = k - n
g_n(k)
\end{array}
$$
where $\lfloor \,
\cdot \,\rfloor$ rounds its argument to the nearest integer towards
minus infinity. We call $g_n$ and $r_n$ index functions; they are
defined on ${\mathbb Z}$ with ranges ${\mathbb Z}$ and $\{1, 2, \ldots,
n\}$, respectively.

If we write
\begin{equation}\label{equ:8-1}
k = j n + i
\end{equation}
with $1 \leq i \leq n$ and $j \in {\mathbb Z}$, then
$$
\begin{array}{rcl} g_n (j n + i) = j &\mbox{and}& r_n(j n
+ i) = i.
\end{array}
$$
%
%
%
%
%
%


\section{ML($n$)BiCG} \label{sec:12-24-1} Analogously to the derivation of
BiCGStab from BiCG, ML($n$)BiCGStab algorithms 
were derived from a BiCG-like
algorithm named ML($n$)BiCG, which was built
upon a band Lanczos process with $n$ left starting vectors and a
single right starting vector. In this section, we present the
ML($n$)BiCG algorithm from \cite{yeungchan}.

Consider the solution of
(\ref{equ:7-6-1}).
Throughout the
paper we do not assume the coefficient matrix ${\bf A}$ is nonsingular.
In \cite{ye2012}, we proved that ML($n$)BiCG/ML($n$)BiCGStab can solve a singular system almost surely provided that the underlining Krylov subspace contains a solution of (\ref{equ:7-6-1}).


Let be given $n$ vectors ${\bf q}_1, \ldots, {\bf
q}_n \in {\mathbb C}^N$, which we call {\it left starting vectors} or {\it shadow vectors}. Define
\begin{equation}
\begin{array}{lll}
 {\bf p}_{k} = \left( {\bf A}^H \right)^{g_n(k)}
{\bf q}_{r_n(k)}, & & k \in {\mathbb N}.
\end{array}
\label{equ:7-9-5}
\end{equation}
The following algorithm for the solution of
(\ref{equ:7-6-1}) is from
\cite{yeungchan}.\\

\begin{algorithm}{\bf ML($n$)BiCG}
\label{alg:1} \vspace{.2cm}
\begin{tabbing}
x\=xxx\= xxx\=xxx\=xxx\=xxx\=xxx\kill \>1. \> Choose an initial
guess $\widehat{\bf x}_0$ and $n$ vectors ${\bf q}_1, {\bf q}_2,
\ldots, {\bf q}_n$. \\
\>2. \>  Compute $\widehat{\bf r}_0 = {\bf b} - {\bf A} \widehat{\bf
x}_0$ and
set ${\bf p}_1 = {\bf q}_1$, $\widehat{\bf g}_0 = \widehat{\bf r}_0$. \\
\>3. \>For $k = 1, 2, \ldots$, until convergence: \\
\>4. \>\>$\alpha_k = {\bf p}_k^H \widehat{\bf r}_{k-1} / {\bf p}_k^H
{\bf A} \widehat{\bf g}_{k-1}$; \\
\>5. \>\>$\widehat{\bf x}_k = \widehat{\bf x}_{k-1} + \alpha_k \widehat{\bf g}_{k-1}$; \\
\>6. \>\> $\widehat{\bf r}_k = \widehat{\bf r}_{k-1} - \alpha_k {\bf A} \widehat{\bf g}_{k-1}$; \\
\>7. \>\>For $s = \max (k - n, 0), \ldots, k - 1$ \\
\>8. \>\>\>$\beta^{(k)}_{s} = - {\bf p}^H_{s+1} {\bf A}
 \left(\widehat{\bf r}_k + \sum_{t = \max (k - n, 0) }^{s-1} \beta^{(k)}_t \widehat{\bf g}_t \right)
\big/{\bf p}^H_{s+1} {\bf A} \widehat{\bf g}_s$; \\
\>9. \>\>End \\
\>10.\>\> $\widehat{\bf g}_k = \widehat{\bf r}_k + \sum_{s = \max (k
- n, 0) }^{k -1} \beta_{s}^{(k)}
            \widehat{\bf g}_{s}$; \\
\>11. \>\>Compute ${\bf p}_{k+1}$ according to (\ref{equ:7-9-5}) \\
\>12. \> End
\end{tabbing}
\end{algorithm}
\vspace{.2cm}

This ML($n$)BiCG algorithm is a variation of the classical BiCG
algorithm with the left-hand side (shadow) Krylov subspace of BiCG being
replaced by the block Krylov subspace
$$\begin{array}{rl}
{\cal B}_k
&
\equiv \mbox{the space spanned by the first } k \mbox{ columns of } [{\bf Q}, {\bf A}^H {\bf Q}, ({\bf A}^H)^2 {\bf Q}, \ldots ]\\
& = span \{{\bf p}_1, {\bf p}_2,
\ldots, {\bf p}_k\}\\
& = \sum_{i = 1}^{r_n(k)}{\cal K}_{g_n(k)+1}({\bf A}^H, {\bf q}_i) +
\sum_{i = r_n(k)+1}^n {\cal K}_{g_n(k)}({\bf A}^H, {\bf q}_i)
\end{array}
$$
where ${\bf Q} 
\equiv [{\bf q}_1, {\bf q}_2, \ldots, {\bf q}_n]$, ${\cal K}_0({\bf M}, {\bf v}) = \{{\bf 0}\}$ and
$$
{\cal K}_t({\bf M}, {\bf v}) 
\equiv span \{{\bf v}, {\bf M} {\bf v},
\ldots, {\bf M}^{t-1} {\bf v}\}
$$
for ${\bf M} \in {\mathbb C}^{N \times N}, {\bf v} \in {\mathbb C}^N$ and
$t \in {\mathbb N}$. Moreover, in this ML($n$)BiCG, the basis used for
${\cal B}_k$ is not chosen to be bi-orthogonal, but simply the set
$\{{\bf p}_1, {\bf p}_2, \ldots, {\bf p}_k\}$. Therefore,
it can be viewed as a generalization of a
one-sided Lanczos algorithm (see \cite{gut10, saad82}).

It can be shown that the quantities of ML($n$)BiCG satisfy the properties (see \cite{ye2012})
\begin{enumerate}
\item[(a)]
$\widehat{\bf x}_k \in \widehat{\bf x}_0 + {\cal K}_k ({\bf A},
\widehat{\bf r}_0)
$, $\widehat{\bf r}_k \in
\widehat{\bf r}_0 + {\bf A} {\cal K}_k ({\bf A}, \widehat{\bf r}_0)
$.
\item[(b)] $
\widehat{\bf r}_k \perp 
span \{ {\bf p}_1, {\bf p}_2, \ldots, {\bf p}_k\}
$
and $\widehat{\bf
r}_k \not\perp {\bf p}_{k+1}$.
\item[(c)] $
{\bf A} \widehat{\bf g}_k \perp 
span \{ {\bf p}_1, {\bf p}_2, \ldots, {\bf p}_k\}
$ and ${\bf A}
\widehat{\bf g}_k \not\perp {\bf p}_{k+1}$.
\end{enumerate}

\section{ML($n$)BiCGStabt}
 \label{sec:12-24-2} The derivation of a ML($n$)BiCGStab algorithm from ML($n$)BiCG essentially is a Sonneveld-van der Vorst-Lanczos procedure. The
central idea
of this procedure is the
remarkable
observation: inner products ${\bf p}^H \widehat{\bf r}$ and ${\bf
p}^H {\bf A} \widehat{\bf g}$ in BiCG can be replaced by inner
products of the forms ${\bf q}^H \psi({\bf A}) \widehat{\bf r}$ and
${\bf q}^H {\bf A} \psi({\bf A}) \widehat{\bf g}$ respectively, where $\psi$ is
an arbitrary polynomial with some suitable degree. This observation can also applied to ML($n$)BiCG because of properties (b) and (c) stated in  \S\ref{sec:12-24-1}.

\subsection{Algorithm}
In \cite{ye2012}, Yeung presented two ML($n$)BiCGStab algorithms, labeled as Algorithms 4.1 and 5.1 respectively. Let $\phi_k$ be the polynomial of degree $k$, recursively defined by
$$
\phi_k(\lambda) = \left\{ \begin{array}{lcl} 1& & \mbox{if } k = 0\\
(1-\omega_k \lambda) \phi_{k-1}(\lambda)& & \mbox{if } k > 0
\end{array} \right.
$$
where $\omega_k$ is a free parameter. Then the
quantities in Algorithm 4.1 are defined by
\begin{equation} \label{equ:7-15-2}
\begin{array}{lcl}
{\bf r}_k = \phi_{g_n(k)+1}({\bf A}) \,\widehat{\bf r}_k,& &{\bf
u}_k = \phi_{g_n(k)} ({\bf A}) \,\widehat{\bf r}_k,\\
{\bf g}_k = \phi_{g_n(k)+1}({\bf A}) \,\widehat{\bf g}_k, & & {\bf
d}_k = -\omega_{g_n(k)+1}{\bf A} \phi_{g_n(k)}({\bf A})\, \widehat{\bf
g}_k,\\
{\bf w}_k = {\bf A} {\bf g}_k
\end{array}
\end{equation}
for $k > 0$, and those in Algorithm 5.1
defined as
\begin{equation}\label{equ:7-15-3}
\begin{array}{lcl}
{\bf r}_k = \phi_{g_n(k+1)}({\bf A}) \,\widehat{\bf r}_k, & & {\bf
g}_k = \phi_{g_n(k+1)}({\bf A}) \widehat{\bf g}_k,\\
 {\bf u}_k = \phi_{g_n(k)}({\bf A}) \widehat{\bf r}_k, & & {\bf w}_k = {\bf A} {\bf g}_k
\end{array}
\end{equation}
for $k > 0$. When $k = 0$, both algorithms set
$$
\begin{array}{lcl}
{\bf r}_0 = \widehat{\bf r}_0 & \mbox{and} & {\bf g}_0 =
\widehat{\bf g}_0.
\end{array}
$$
Here
${\bf r}_k$ is the residual
of the $k$th
approximate solution ${\bf x}_k$. Numerical experiments in \cite{ye2012} indicated that the ${\bf r}_k$ computed by Algorithm 4.1 is generally closer to the true residual ${\bf b} - {\bf A} {\bf x}_k$ than the ${\bf r}_k$ computed by Algorithm 5.1 is. A close examination of the algorithms can explain this difference in stability.

In both algorithms, the ${\bf x}_k$ and ${\bf r}_k$ are updated by the recursive relations
$$
\begin{array}{rcl}
{\bf x}_k = {\bf x}_{k-1} + \alpha_k {\bf g}_{k-1}, & &{\bf r}_k = {\bf r}_{k-1} - \alpha_k {\bf w}_{k-1}
\end{array}
$$
in most $k$-iterations, where $\alpha_k$ is a scalar. The true residual of the computed ${\bf x}_k$ is therefore
\begin{equation}\label{equ:7-15-10}
{\bf b} - {\bf A} {\bf x}_k = {\bf b} - {\bf A} ({\bf x}_{k-1} + \alpha_k {\bf g}_{k-1}) = ({\bf b} - {\bf A} {\bf x}_{k-1}) - \alpha_k {\bf A}{\bf g}_{k-1}.
\end{equation}
In Algorithm 4.1, ${\bf w}_k$ is updated by ${\bf w}_k = {\bf A}{\bf g}_k$ (as it is defined in (\ref{equ:7-15-2})) 
for all $k$.
In Algorithm 5.1, however, ${\bf w}_k$ is updated by ${\bf w}_k = {\bf A}{\bf g}_k$ only when $r_n(k) = n$.
In other words,
the update for ${\bf r}_k$ in Algorithm 4.1 follows closer to (\ref{equ:7-15-10})\footnote{There is a similar comment on BiCGStab($l$) \cite[p.27]{SF} when compared to BiCGStab2\cite{gut}.}.
As a result,
the
residual ${\bf r}_k$ computed by Algorithm 4.1 is generally closer to the true residual (\ref{equ:7-15-10}) than the ${\bf r}_k$ computed by Algorithm 5.1 is.
Because of the observation, we expect that Algorithm 5.1 should be as stable as Algorithm 4.1 if we could modify the algorithm so that its ${\bf w}_k$ were updated by ${\bf w}_k = {\bf A}{\bf g}_k$ in all the iterations --- this is the goal that we develop ML($n$)BiCGStabt.

The derivation of Algorithm 5.1 in \cite{ye2012} was divided into several stages, starting from ML($n$)BiCG. The following is a copy of its Derivation Stage \#8 which is a list of equations that the quantities in (\ref{equ:7-15-3}) satisfy.\\

{\bf Derivation Stage \#8 in \cite{ye2012}.} \vspace{.2cm}
\begin{tabbing}
x\=xxx\= xxx\=xxx\=xxx\=xxx\=xxx\=xxx\=xxx\=xxx\kill
\>1. \>For $k = 1, 2, \ldots$, until convergence: \\
\>2. \>\>$\alpha_k = {\bf q}_{r_n(k)}^H {\bf r}_{k-1} / {\bf
q}_{r_n(k)}^H
 {\bf w}_{k-1}$;\\
\>3.\>\> If $r_n(k) < n$\\
\>4. \>\>\> ${\bf r}_k = {\bf r}_{k-1} - \alpha_k {\bf w}_{k-1}$;
\\
\>5. \>\>\>For $s = \max (k - n, 0), \ldots, g_n(k)n - 1$ \\
\>6. \>\>\>\>$\beta^{(k)}_{s} = {\bf q}^H_{r_n(s+1)}
 \left(
  {\bf r}_k - \,\omega_{g_n(k+1)} \sum_{t = \max (k - n, 0)
}^{s-1} \beta^{(k)}_t {\bf w}_t \right) \big/\omega_{g_n(k+1)}
 {\bf q}^H_{r_n(s+1)} {\bf w}_s$; \\
\>7. \>\>\>End
\\
\>8. \>\>\>For $s = g_n(k)n, \ldots, k - 1$ \\
\>9. \>\>\>\>$\beta^{(k)}_{s} = - {\bf q}^H_{r_n(s+1)} \left({\bf A}
 {\bf
r}_k + \sum_{t = \max (k - n, 0) }^{g_n(k)n-1} \beta^{(k)}_t
({\bf I} - \omega_{g_n(k+1)} {\bf A})
 {\bf w}_t
  \right.$\\
  \>\>\>\>\>\>\> $\left. + \sum_{t = g_n(k)n }^{s-1} \beta^{(k)}_t
 {\bf w}_t
 \right)
\big/ {\bf q}^H_{r_n(s+1)} {\bf w}_s$; \\
\>10. \>\>\>End
\\
\>11.\>\>\> ${\bf g}_k = {\bf r}_k - \omega_{g_n(k+1)} \sum_{s = \max
(k - n, 0) }^{g_n(k)n -1} \beta_{s}^{(k)} {\bf w}_{s} + \sum_{s =
\max (k - n, 0) }^{g_n(k)n -1} \beta_{s}^{(k)} {\bf g}_{s} + \sum_{s
= g_n(k)n }^{k -1} \beta_{s}^{(k)} {\bf g}_{s}
$;\\
\>12.\>\>Else\\
\>13. \>\>\> $ {\bf u}_k = {\bf r}_{k-1} - \alpha_k {\bf w}_{k-1}$;
\\
\>14. \>\>\>${\bf r}_k = ({\bf I} - \omega_{g_n(k+1)} {\bf A})
{\bf u}_k$;\\
\>15. \>\>\>For $s = g_n(k)n, \ldots, k - 1$ \\
\>16. \>\>\>\>$\beta^{(k)}_{s} = {\bf q}^H_{r_n(s+1)} \left( {\bf
r}_k - \omega_{g_n(k+1)} \sum_{t = g_n(k)n }^{s-1} \beta^{(k)}_t
  {\bf w}_t
 \right)
\big/\omega_{g_n(k+1)}
 {\bf q}^H_{r_n(s+1)} {\bf w}_s$;\\
\>17. \>\>\>End \\
\>18.\>\>\> ${\bf g}_k = {\bf r}_k - \omega_{g_n(k+1)} \sum_{s =
g_n(k)n }^{k -1} \beta_{s}^{(k)} {\bf w}_{s} + \sum_{s = g_n(k)n
}^{k -1} \beta_{s}^{(k)} {\bf g}_{s}
$; \\
\>19.\>\>End\\
\>20. \> End
\end{tabbing}

\vspace{.2cm}

According to (\ref{equ:7-15-3}), the equation in Line 9 can be rewritten as
$$\begin{array}{rcl}
\beta^{(k)}_{s} &=& - {\bf q}^H_{r_n(s+1)} {\bf A} \left(
 {\bf
r}_k + \sum_{t = \max (k - n, 0) }^{g_n(k)n-1} \beta^{(k)}_t
({\bf I} - \omega_{g_n(k+1)} {\bf A})
 {\bf g}_t \right.\\
 & & \left.
   + \sum_{t = g_n(k)n }^{s-1} \beta^{(k)}_t
 {\bf g}_t
 \right)
\big/ {\bf q}^H_{r_n(s+1)} {\bf w}_s\\
&=& - {\bf q}^H_{r_n(s+1)} {\bf A} \left(
 {\bf
r}_k + \sum_{t = \max (k - n, 0) }^{g_n(k)n-1} \beta^{(k)}_t
({\bf g}_t - \omega_{g_n(k+1)} {\bf w}_t)\right.\\
& & \left.
   + \sum_{t = g_n(k)n }^{s-1} \beta^{(k)}_t
 {\bf g}_t
 \right)
\big/ {\bf q}^H_{r_n(s+1)} {\bf w}_s.
\end{array}
$$
It is because of the $\bf A$ before the parentheses, we can not update ${\bf w}_k$ by ${\bf w}_k = {\bf A}{\bf g}_k$ in Algorithm 5.1 while keeping the average number of matrix-vector multiplications as low as $1 + 1/n$ per iteration. If, however, the vector ${\bf f}_{r_n(s+1)} \equiv {\bf A}^H {\bf q}_{r_n(s+1)}$ is available, then Line 9 will become
\begin{equation}\label{equ:7-16-1}
\begin{array}{rcl}
\beta^{(k)}_{s}
&=& - {\bf f}^H_{r_n(s+1)} \left(
 {\bf
r}_k + \sum_{t = \max (k - n, 0) }^{g_n(k)n-1} \beta^{(k)}_t
({\bf g}_t - \omega_{g_n(k+1)} {\bf w}_t)\right.\\
& &
 \left.  + \sum_{t = g_n(k)n }^{s-1} \beta^{(k)}_t
 {\bf g}_t
 \right)
\big/ {\bf q}^H_{r_n(s+1)} {\bf w}_s
\end{array}
\end{equation}
and the troubling $\bf A$ is gone. It is the observation that leads to the ML($n$)BiCGStabt algorithm.

Replace Line 9 in Derivation Stage \#8 with (\ref{equ:7-16-1}) and suppose
$$
{\bf F} \equiv {\bf A}^H {\bf Q} = [{\bf A}^H {\bf q}_1, {\bf A}^H {\bf q}_2, \ldots, {\bf A}^H {\bf q}_n]
$$
is available.
Recalling that ${\bf r}_k$ is the residual of ${\bf x}_k$,
to be consistent
with Lines 4, 13 and 14, we update the approximate solution ${\bf x}_k$
as
\begin{equation}\label{equ:7-16-2}
\begin{array}{rl}
{\bf x}_k = & \left\{ \begin{array}{lcl} {\bf x}_{k-1} + \alpha_{k}
 {\bf g}_{k-1}, & &\mbox{if } r_n(k) < n \\
 \omega_{g_n(k+1)} {\bf u}_k + {\bf x}_{k-1} + \alpha_k {\bf g}_{k-1}, & & \mbox{if } r_n(k) = n.
\end{array} \right.
\end{array}
\end{equation}
Now adding (\ref{equ:7-16-2}) and ${\bf w}_k = {\bf A} {\bf g}_k$ to the derivation stage, then simplifying the operations
appropriately, we arrive at the following algorithm. The free
parameter $\omega_{g_n(k+1)}$ is chosen to minimize the $2$-norm of
${\bf r}_k$.\\

\begin{algorithm}{\bf ML($n$)BiCGStabt without preconditioning 
} \label{alg:13}
\vspace{.2cm}
\begin{tabbing}
x\=xxx\= xxx\=xxx\=xxx\=xxx\=xxx\=xxx\=xxx\=xxx\kill
\>1. \> Choose an initial guess ${\bf
x}_0$ and $n$ vectors ${\bf q}_1, {\bf q}_2, \ldots,
{\bf q}_n$. \\
\>2. \> Compute $[{\bf f}_1, \ldots, {\bf f}_{n-1}] = {\bf A}^H [{\bf q}_1, \ldots, {\bf q}_{n-1}]$. \\
\>3. \>  Compute ${\bf r}_0 = {\bf b} - {\bf A} {\bf x}_0$ and ${\bf
g}_0 = {\bf r}_0,\,\, {\bf w}_0 = {\bf A}
{\bf g}_0,\,\, c_0 = {\bf q}^H_{1} {\bf w}_0$, $\omega_0 = 1$. \\
\>4. \>For $k = 1, 2, \ldots$, until convergence: \\
\>5. \>\>$\alpha_k = {\bf q}_{r_n(k)}^H {\bf r}_{k-1} / c_{k-1}$;\\
\>6.\>\> If $r_n(k) < n$\\
\>7. \>\>\> ${\bf x}_k = {\bf x}_{k-1} + \alpha_k {\bf g}_{k-1}$;
\,
${\bf r}_k = {\bf r}_{k-1} - \alpha_k {\bf w}_{k-1}$;
\\
\>8. \>\>\> ${\bf z}_w = {\bf r}_{k}; \,\, {\bf g}_k = {\bf 0}$;
\\
\>9. \>\>\>For $s = \max (k - n, 0), \ldots, g_n(k)n - 1$ \\
\>10. \>\>\>\>$\tilde{\beta}^{(k)}_{s} = - {\bf q}^H_{r_n(s+1)}
 {\bf z}_w \big/
 c_s$; \,\,\,\,\,\, \% $\tilde{\beta}^{(k)}_{s} = -\omega_{g_n(k+1)} \beta^{(k)}_{s}$\\
 \>11. \>\>\>\>${\bf z}_w = {\bf z}_w + \tilde{\beta}^{(k)}_{s} {\bf w}_s$;\\
 \>12. \>\>\>\>${\bf g}_k = {\bf g}_k + \tilde{\beta}^{(k)}_{s} {\bf g}_s$; \\
\>13. \>\>\>End
\\
\>14. \>\>\>${\bf g}_k = {\bf z}_w - \frac{1}{\omega_{g_n(k+1)}} {\bf g}_k$;
\\
\>15. \>\>\>For $s = g_n(k)n, \ldots, k - 1$ \\
\>16. \>\>\>\>$\beta^{(k)}_{s} = - {\bf f}_{r_n(s+1)}^H {\bf g}_k
\big/ c_s$; \\
\>17. \>\>\>\>${\bf g}_k = {\bf g}_k + \beta^{(k)}_{s} {\bf g}_s$; \\
\>18. \>\>\>End
\\
\>19.\>\>Else\\
\>20. \>\>\> $ {\bf x}_k = {\bf x}_{k-1} + \alpha_k {\bf g}_{k-1}$;
\\
\>21. \>\>\>
$ {\bf u}_k = {\bf r}_{k-1} - \alpha_k {\bf w}_{k-1}$;
\\
\>22. \>\>\> $\omega_{g_n(k+1)} =  ({\bf A} {\bf u}_k)^H {\bf u}_k / \|{\bf A}{\bf u}_k \|_2^2$;\\
\>23.\>\>\>${\bf x}_k = {\bf x}_k  +\omega_{g_n(k+1)} {\bf u}_k$;\,
${\bf r}_k = -\omega_{g_n(k+1)} {\bf A}{\bf u}_k +
{\bf u}_k$; \\
\>24. \>\>\> ${\bf z}_w = {\bf r}_{k}; \,\, {\bf g}_k = {\bf 0}$;
\\
\>25. \>\>\>For $s = g_n(k)n, \ldots, k - 1$ \\
\>26. \>\>\>\>$\tilde{\beta}^{(k)}_{s} = - {\bf q}^H_{r_n(s+1)} {\bf z}_w
\big/
 c_s$; \,\,\,\,\,\, \% $\tilde{\beta}^{(k)}_{s} = -\omega_{g_n(k+1)} \beta^{(k)}_{s}$\\
 \>27. \>\>\>\>${\bf z}_w = {\bf z}_w + \tilde{\beta}^{(k)}_{s} {\bf w}_s$;\\
 \>28. \>\>\>\>${\bf g}_k = {\bf g}_k + \tilde{\beta}^{(k)}_{s} {\bf g}_s$;\\
\>29. \>\>\>End \\
\>30.\>\>\> ${\bf g}_k = {\bf z}_w  - \frac{1}{\omega_{g_n(k+1)}} {\bf g}_{k}
$; \\
\>31.\>\>End\\
\>32.\>\> ${\bf w}_k = {\bf A} {\bf g}_k
$; \,
$c_k = {\bf q}_{r_n(k+1)}^H {\bf w}_k$;\\
\>33. \> End
\end{tabbing}
\end{algorithm}\vspace{.2cm}

Line 32 indicates that ${\bf w}_k$ is computed by ${\bf w}_k = {\bf A} {\bf g}_k$ for all $k$-iterations. Therefore the updates
for ${\bf x}_k$ and ${\bf r}_k$ in the above Algorithm \ref{alg:13} are
\begin{equation}\label{equ:9-5-2}
\begin{array}{ccc}
{\bf x}_k = {\bf x}_{k-1} + \alpha_k {\bf g}_{k-1}, & & {\bf r}_k = {\bf r}_{k-1} - \alpha_k {\bf A} {\bf g}_{k-1}
\end{array}
\end{equation}
which meets the goal that we set right before Derivation Stage \#8 on improving the stability of Algorithm 5.1 in \cite{ye2012}.
The stability of updates of the type (\ref{equ:9-5-2}) has been studied in detail by Neumaier\cite{neum} and Sleijpen and van der Vorst\cite{Svan}.

We remark that (i) Algorithm \ref{alg:13} does not compute ${\bf u}_k$ when
$r_n(k) < n$. In fact, ${\bf u}_k = {\bf r}_k$ when $r_n(k) < n$ from (\ref{equ:7-15-3}); (ii) if the ${\bf
u}_k$ in Line 21 happens to be zero, then the ${\bf x}_k$ in Line 20
will be the exact solution to system (\ref{equ:7-6-1}) and the
algorithm stops there.

Computational and storage cost
based on the preconditioned ML($n$)BiCGStabt (see Algorithm \ref{alg:10-22})
is presented in Table \ref{tab:7-16-1}. Note that we do not need to store both ${\bf A}$ and ${\bf A}^H$ since ${\bf A}^H$ is only used in Line 2. Compared with Algorithm 5.1 in \cite{ye2012}, the computational cost of ML($n$)BiCGStabt is slightly cheaper.

\begin{table}[tbp]
\center
\caption{Average cost per $k$-iteration of the
preconditioned ML($n$)BiCGStabt Algorithm \ref{alg:10-22} and its storage. This table does not count the cost in Lines 1-2 of the algorithm.}
\begin{tabular}{||c|c|c|c||}  \hline
Preconditioning ${\bf M}^{-1} {\bf v}$ & $\displaystyle{1 +
1/n}$  & ${\bf u} \pm {\bf v}, \,\, \alpha {\bf v}$ &
$\displaystyle{1}$ \\ \hline Matvec ${\bf A} {\bf v}$ &
$\displaystyle{1 + 1/n}$  & Saxpy ${\bf u} + \alpha {\bf
v}$ & $\displaystyle{1.5 n + 2.5 + 2/n}$  \\ \hline dot product
$\displaystyle{ {\bf u}^H {\bf v}}$ & $\displaystyle{n + 1+
2/n}$ & Storage &${\bf A} + {\bf M} +$
\\
 & & & $(4 n+4) N +O(n)$\\ \hline
\end{tabular}\label{tab:7-16-1}
\end{table}

Theoretically, it can be guaranteed that an exact breakdown in Algorithm \ref{alg:13} is almost impossible (see \cite{ye2012} for a detailed analysis). The algorithm, however, can encounter a near breakdown in its implementation. The divisors in the algorithm are $c_k, \|{\bf A} {\bf u}_k\|_2$ and $\omega_{g_n(k+1)}$. If $\|{\bf A} {\bf u}_k\|_2 \approx 0$, then ${\bf u}_k \approx {\bf 0}$ and the ${\bf x}_k$ in Line 20 is an approximate solution. If $\omega_{g_n(k+1)} \approx 0$, we can add some small perturbation to it so that it is relatively far from $0$. About $c_k$,
it can be showed that it is a quantity that relates to $\omega_{g_n(k+1)}$ and the ML($n$)BiCG divisor ${\bf p}_{k+1}^H {\bf A} \widehat{\bf g}_k$. The ML($n$)BiCG divisor ${\bf p}_{k+1}^H {\bf A} \widehat{\bf g}_k$ is in turn related to the underlying Lanczos breakdown and the breakdown caused by the non-existence of the $LU$ factorization of the Hessenberg matrix of the recurrence coefficients. But, as indicated in \cite{gut}, in most cases such breakdowns can be overcome by a look-ahead step, see \cite{fgn, gut30, gut31, parlett} and further references cited there. Moreover, for how to avoid a breakdown in a nonsymmetric block Lanczos algorithm, one can consult \cite{loher}.

\subsection{Properties} Since the quantities of ML($n$)BiCGStabt are defined exactly the same as those of Algorithm 5.1 in \cite{ye2012}, ML($n$)BiCGStabt shares the same properties with Algorithm 5.1.

Let $\nu$ be the degree of the minimal polynomial $p_{min}(\lambda;
{\bf A}, {\bf r}_0)$ of ${\bf r}_0$ with respect to
${\bf A}$, namely, the unique monic polynomial $p(\lambda)$ of
minimum degree such that $p({\bf A}) {\bf r}_0 = {\bf 0}$,
and let
$$
{\bf S}_\nu = [{\bf p}_1, {\bf p}_2, \ldots, {\bf p}_\nu]^H {\bf A}
[{\bf r}_0, {\bf A} {\bf r}_0, \ldots,
{\bf A}^{\nu-1} {\bf r}_0]
$$
and
$$
{\bf W}_\nu = [{\bf p}_1, {\bf p}_2, \ldots, {\bf p}_\nu]^H
[{\bf r}_0, {\bf A} {\bf r}_0, \ldots, {\bf
A}^{\nu-1} {\bf r}_0].
$$
Denote by ${\bf S}_l$ and ${\bf W}_l$ the $l \times
l$ leading principal submatrices of ${\bf S}_\nu$ and
${\bf W}_\nu$ respectively.
Joubert\cite{joubert1, joubert2} called these matrices the moment matrices. With the notations,
some facts about ML($n$)BiCGStabt (Algorithm \ref{alg:13})
are summarized as follows.\\

\begin{proposition} \cite[Prop.~5.1]{ye2012} In infinite precision arithmetic, if $
\prod_{l =
1}^\nu \det({\bf S}_l) \det({\bf W}_l) \ne 0
$,  $\omega_{g_n(k+1)} \ne 0$ and $
1/\omega_{g_n(k+1)} \not\in \sigma({\bf A})$ for $1 \leq k \leq
\nu-1$, where $\sigma({\bf A})$ is the spectrum of $\bf A$,
then
Algorithm \ref{alg:13} does not break down by zero division for $k =
1, 2, \ldots, \nu$, and the approximate solution ${\bf x}_\nu$ at
step $k = \nu$ is exact to the system (\ref{equ:7-6-1}). Moreover,
the computed quantities satisfy
\begin{enumerate}
\item[(a)]
${\bf x}_k \in {\bf x}_0 + span \{ {\bf r}_0, {\bf A} {\bf r}_0,
\ldots, {\bf A}^{g_n(k+1)+k-1} {\bf r}_0 \}$ and ${\bf r}_k = {\bf
b} - {\bf A} {\bf x}_k \in {\bf r}_0 + span \{ {\bf A} {\bf r}_0,
{\bf A}^2 {\bf r}_0, \ldots, {\bf A}^{g_n(k+1)+k} {\bf r}_0 \}$ for
$1 \leq k \leq \nu-1$.

\item[(b)] ${\bf r}_k \ne {\bf 0}$ for $1 \leq k \leq \nu -1$; ${\bf
r}_\nu = {\bf 0}$.

\item[(c)] ${\bf r}_k \perp span \{ {\bf q}_1, {\bf q}_2, \ldots,
{\bf q}_{r_n(k)}\}$ and ${\bf r}_k \not\perp {\bf q}_{r_n(k)+1}$ for
$1 \leq k \leq \nu -1$ with $r_n(k) < n$; ${\bf r}_k \not\perp {\bf
q}_{1}$ for $1 \leq k \leq \nu -1$ with $r_n(k) = n$.

\item[(d)] ${\bf u}_k \perp span \{ {\bf q}_1, {\bf q}_2, \ldots,
{\bf q}_{n}\}$ for $1 \leq k \leq \nu$ with $r_n(k) = n$.


\item[(e)] ${\bf A}{\bf g}_k \perp span \{ {\bf q}_1, {\bf q}_2, \ldots,
{\bf q}_{r_n(k)}\}$ and ${\bf A}{\bf g}_k
 \not\perp {\bf q}_{r_n(k)+1}$ for
$1 \leq k \leq \nu -1$ with $r_n(k) < n$; ${\bf A} {\bf g}_k
\not\perp {\bf q}_{1}$ for $1 \leq k \leq \nu -1$ with $r_n(k) = n$.
\end{enumerate}
\end{proposition}\vspace{.2cm}

In \S6.2 of \cite{ye2012}, relations of Algorithm 5.1 in \cite{ye2012} to some existing methods were presented. The same arguments applied to ML($n$)BiCGStabt imply that
\begin{enumerate}
\item[(a)] ML($n$)BiCGStabt is a FOM algorithm, but involving ${\bf A}^H$ in its implementation, if we set $n \geq \nu$ and ${\bf q}_k = {\bf r}_{k-1}$.
    \item[(b)] ML($n$)BiCGStabt is a BiCGStab algorithm if we set $n = 1$.
    \item[(c)] ML($n$)BiCGStabt is a IDR($s$) algorithm with $s = n$, but involving ${\bf A}^H$ in its implementation.
\end{enumerate}

\section{Numerical experiments}\label{sec:12-19}
A preconditioned ML($n$)BiCGStabt algorithm can be obtained by
applying Algorithm \ref{alg:13} to
the system
$$
{\bf A} {\bf M}^{-1} {\bf y} = {\bf b}
$$
where ${\bf M}$ is nonsingular, then recovering ${\bf x}$ through
${\bf x} = {\bf M}^{-1} {\bf y}$. The resulting algorithm,
Algorithm \ref{alg:10-22}, together
with its Matlab code are presented in \S\ref{sec:apen}. To avoid calling the index
functions $r_n(k)$ and $g_n(k)$ every $k$-iteration, we have split
the $k$-loop into a $i$-loop and a $j$-loop where $i, j, k$ are
related by (\ref{equ:8-1}) with $1 \leq i \leq n, 0 \leq j$.
Moreover, we have optimized the operations as much as possible in
the resulting preconditioned algorithm.

We compared ML($n$)BiCGStabt with BiCG,
 BiCGStab and two algorithms of ML($n$)BiCGStab: Algorithms 4.1 and 5.1 in \cite{ye2012}.
All test data were downloaded from The University of Florida Sparse Matrix Collection\footnote{http://www.cise.ufl.edu/research/sparse/matrices/}, and
the computing was done in Matlab Version 7.1 on a Windows XP
machine with a Pentium 4 processor.
In all the experiments, we chose the initial guess ${\bf x}_0 = {\bf 0}$, the
stopping criterion $\| {\bf r}_k \|_2 / \| {\bf b}\|_2 < 10^{-7}$
where ${\bf r}_k$ was the computed residual, and the Sleijpen-van der Vorst minimization control parameter (see \cite{ye2012}) $\kappa = 0$. As for the shadow vectors, we chose
${\bf Q} = [{\bf r}_0, sign(randn(N,n-1))]$.
When a data did not provide a right-hand side, we set ${\bf b} = {\bf A} {\bf e}$ where ${\bf e}$ is the vector of ones.

{\it Example 1.}
We ran all the methods on the
selected
group of matrices
in Table \ref{tab:50}. No preconditioner was used.
The results are summarized in Tables \ref{tab:18}-\ref{tab:8-2-1}. The ``True error'' columns in the tables contain the true relative errors $\|{\bf b} - {\bf A}{\bf x}\|_2/\|{\bf b}\|_2$ where ${\bf x}$ is the computed solution output by an algorithm when it converges. In this experiment, we observe that ML($n$)BiCGStabt and ML($n$)BiCGStab generally outperform BiCG and BiCGStab in terms of computational time.
As an improved version of Algorithm 5.1 in \cite{ye2012}, ML($n$)BiCGStabt has the same stability with Algorithm 4.1 in \cite{ye2012} and is slightly more stable than Algorithm 5.1.

{\it Example 2.} Our experience with the Florida collection has shown that Algorithm 5.1 in \cite{ye2012} is overall a stable algorithm. But still, one can find one or two matrices where it is unstable. Consider the data

\begin{enumerate}
\item {\it e40r0100}, a 2D/3D problem from the Shen group. The coefficient matrix is a $17281$-by-$17281$ real unsymmetric matrix with
$553,562$ nonzero entries.
\item {\it utm5940}, an electromagnetics problem from the TOKAMAK group. The coefficient matrix is a
$5940$-by-$5940$ real unsymmetric matrix with
$83,842$ nonzero entries.
\end{enumerate}

\begin{figure}[htbp]
\centerline{\hbox{ \psfig{figure=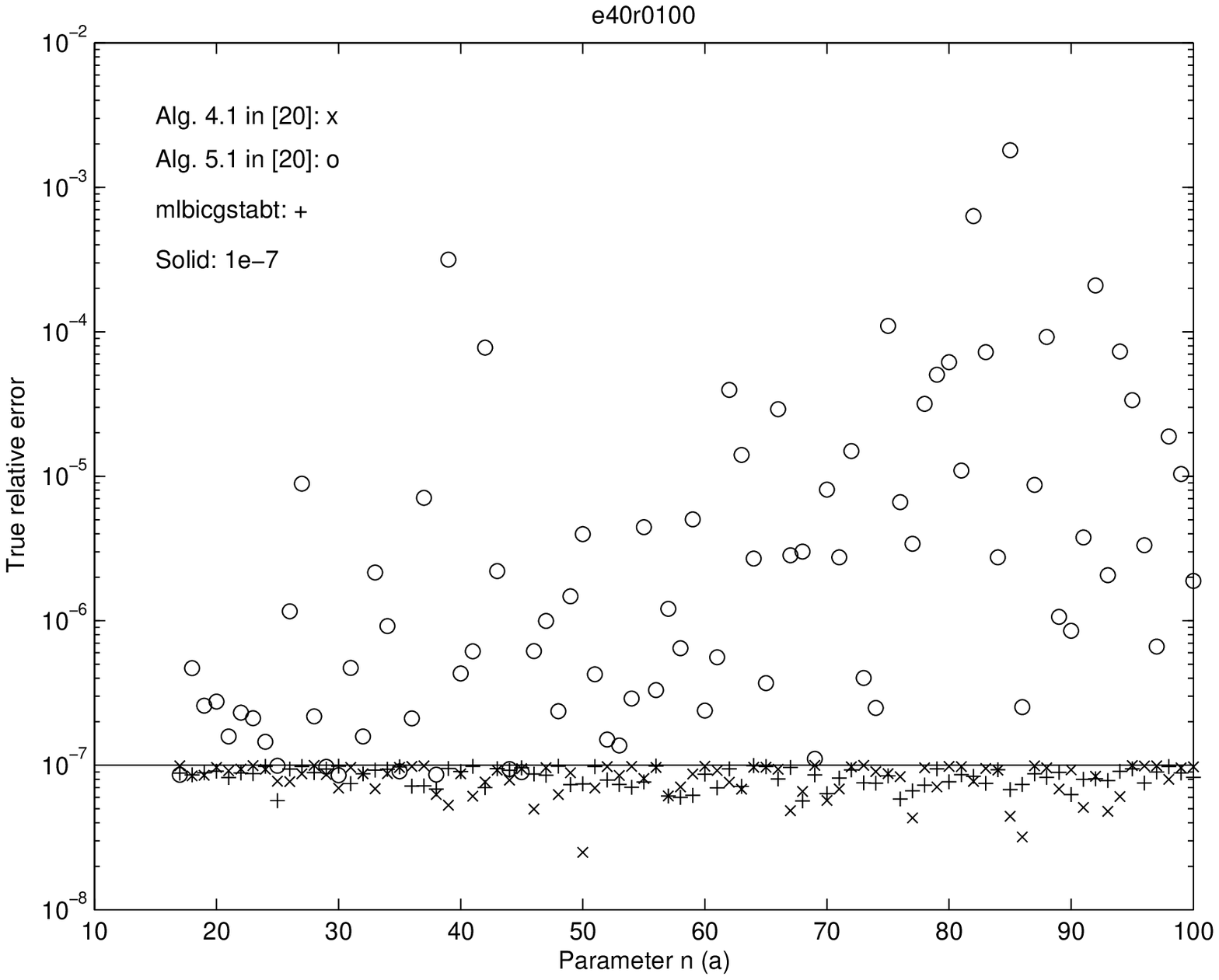,height=6cm }
\psfig{figure=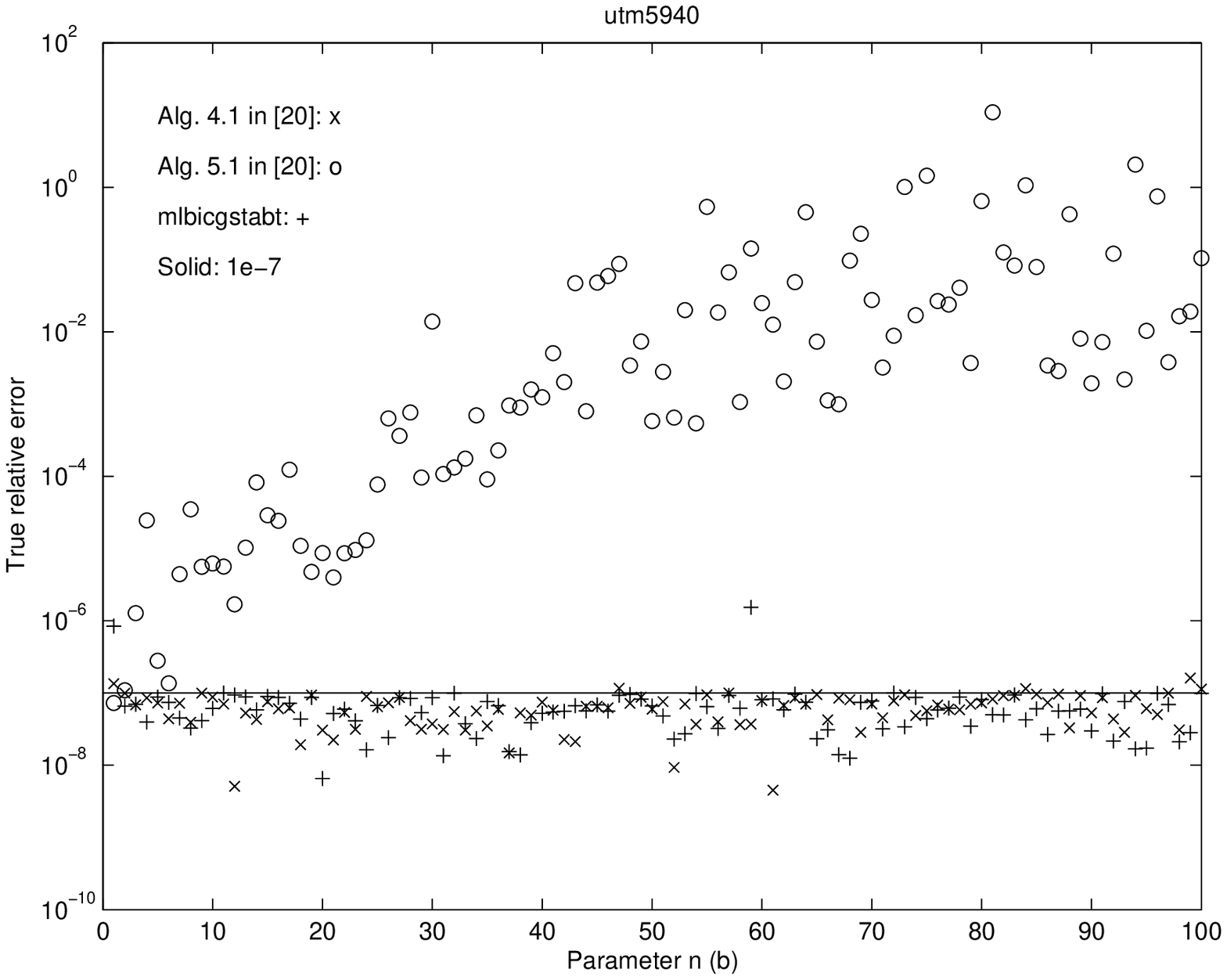,height=6cm}
 }} \caption{Graphs of $E(n)$ against $n$. (a) e40r0100. Since all the three algorithms do not converge when $1 \leq n \leq 16$, we only plot the graphs over the range $17 \leq n \leq 100$; (b) utm5940.
}
\label{fig:11-24-1}
\end{figure}

In this experiment,
ILU preconditioners generated by the Matlab command
[L,U,P] = luinc(A, 1e-3) were used.
For the ease of presentation, we introduce the true relative error
function
$E(n)
\equiv \| {\bf b} - {\bf
A} {\bf x}\|_2 / \|{\bf b}\|_2$ where ${\bf x}$ is the computed solution output by a
ML($n$)BiCGStab algorithm when it converges.
The graphs of $E(n)$ are plotted in Figure \ref{fig:11-24-1}.
It can be seen that
the computed relative errors $\|{\bf r}_k\|_2/\|{\bf b}\|_2$ by Algorithm 5.1 significantly
diverge from their exact counterparts.
By contrast, however, the computed $\|{\bf r}_k\|_2 / \| {\bf
b}\|_2$ by ML($n$)BiCGStabt and Algorithm 4.1 in \cite{ye2012} well approximate their
corresponding true relative errors.
In this experiment, the improvement on stability of ML($n$)BiCGStabt over Algorithm 5.1 is significant.


\begin{table}
\centering
\caption{
A group of data selected from the Florida collection.
Data \#13 contains multiple right-hand sides, and we selected the $45$th in our experiments.
}
\footnotesize{
{\newcommand{\q}[1]{\mc{1}{|l|}{\small\tt #1}}
\noindent
\begin{tabular}{*{5}{|c}|} \hline
\mcc{No.}
&\mcc{Matrix name}&
\mcc{Group name} & \mcc{Size}&\mcc{Nonzeros}
\\ \hline\hline
\mcc{1}& rdb5000&
   \mcc{Bai} & $5,000$ & $29,600$
\\ \hline
\mcc{2}& \mcc{sherman3}&
   \mcc{HB} & $5,005$& $20,033$
\\ \hline
\mcc{3}& \mcc{olm5000}&
   \mcc{Bai} & $5,000$& $19,996$
\\ \hline
\mcc{4}& \mcc{cavity19}&
   \mcc{Drivcav}  & $4,562$& $131,735$
\\ \hline
\mcc{5} & \mcc{tols4000}&
   \mcc{Bai}
   & $4,000$& $8,784$
\\ \hline
\mcc{6} & \mcc{ex31}&
   \mcc{Fidap}  & $3,909$& $91,223$
\\ \hline
\mcc{7}& \mcc{sherman5} &
   \mcc{HB}  & $3,312$& $20,793$
\\ \hline
\mcc{8}& \mcc{raefsky2}&
   \mcc{Simon}   & $3,242$& $293,551$
\\ \hline
\mcc{9}& \mcc{garon1}&
   \mcc{Garon}  & $3,175$& $84,723$
   \\ \hline
\mcc{10}& \mcc{utm5940}&
   \mcc{Tokamak}
   & $5,940$ & $83,842$
\\ \hline
   \mcc{11}& \mcc{Chebyshev3}&
   \mcc{Muite}
   & $4,101$ & $36,879$
\\ \hline
\mcc{12}& \mcc{pores\_2}&
   \mcc{HB}   & $1,224$ & $9,613$
\\ \hline
\mcc{13}& \mcc{tsopf\_rs\_b162\_c1}&
   \mcc{Tsopf}
   & $5,374$ & $205,399$
\\ \hline
\mcc{14}& \mcc{rw5151}&
   \mcc{Bai}
   & $5,151$ & $20,199$
\\ \hline
\mcc{15}& \mcc{circuit\_2}&
   \mcc{Bomhof}
   & $4,510$ & $21,199$
\\ \hline
\mcc{16}& \mcc{viscoplastic1}&
   \mcc{Quaglino}
   & $4,326$ & $61,166$
\\ \hline
\mcc{17}& \mcc{heart1}&
   \mcc{Norris}
   & $3,557$ & $1,385,317$
\\ \hline
\mcc{18}& \mcc{cage9}&
   \mcc{vanHeukelum}
   & $3,534$ & $41,594$
\\ \hline
\mcc{19}& \mcc{thermal}&
   \mcc{Brunetiere}
   & $3,456$ & $66,528$
\\ \hline
\mcc{20}& \mcc{raefsky6}&
   \mcc{Simon}
   & $3,402$ & $130,371$
\\ \hline
\end{tabular}}
}
\label{tab:50}
\end{table}

\begin{table}
\centering
\caption{Experimental results run on the data in Table \ref{tab:50}.
``$-$'' means no convergence within $10N$
iterations.
}
\footnotesize{
{\newcommand{\q}[1]{\mc{1}{|l|}{\small\tt #1}}
\noindent
\begin{tabular}{*{8}{|c}|} \hline
&
\mc{3}{|c|}{BiCG}& \mc{3}{|c|}{BiCGStab}
\\ \hline
No.
& Iter & Time (s) & True error & Iter & Time (s) & True error
\\ \hline\hline
 1
   &$190$  &$0.3823$ & $8.6940 \times 10^{-8}$ & $243$ &$0.2611$
   & $2.4796 \times 10^{-8}$
\\ \hline
2&$-$  &$-$ & $-$ & $6,056$ &$5.2097$
   & $8.6957 \times 10^{-8}$
\\ \hline
3&$5,938$  &$9.2454$ & $7.7274 \times 10^{-8}$ & $-$ &$-$
   & $-$
\\ \hline
4&$24,125$  &$192.3187$ & $6.4603 \times 10^{-8}$ & $-$ &$-$
   & $-$
\\ \hline
5&$-$  &$-$ & $-$ & $-$ &$-$
   & $-$
\\ \hline
6&$3,351$  &$18.6821$ & $5.2792 \times 10^{-8}$ & $4,030$ &$10.2566$
   & $6.7891 \times 10^{-8}$
\\ \hline
7&$1,785$  &$2.6826$ & $8.1858 \times 10^{-8}$ & $3,148$ &$2.4763$
   & $7.7845 \times 10^{-8}$
\\ \hline
8 &$366$  &$6.5125$ & $2.6274 \times 10^{-8}$ & $-$ &$-$
   & $-$
\\ \hline
9&$3,131$  &$15.7303$ & $9.9857 \times 10^{-8}$ & $-$ &$-$
   & $-$
\\ \hline
10&$10,900$  &$57.4226$ & $9.0982 \times 10^{-8}$ & $-$ &$-$
   & $-$
\\ \hline
11&$-$  &$-$ & $-$ & $-$ &$-$
   & $-$
\\ \hline
12&$9,371$  &$6.2491$ & $9.9705 \times 10^{-8}$ & $-$ &$-$
   & $-$
\\ \hline
13&$-$  &$-$ & $-$ & $-$ &$-$
   & $-$
\\ \hline
14&$18$  &$0.0294$ & $7.8645 \times 10^{-8}$ & $-$ &$-$
   & $-$
\\ \hline
15&$-$  &$-$ & $-$ & $334$ &$0.2818$
   & $3.5640 \times 10^{-8}$
\\ \hline
16&$375$  &$1.4655$ & $9.9070 \times 10^{-8}$ & $-$ &$-$
   & $-$
\\ \hline
17&$-$  &$-$ & $-$ & $-$ &$-$
   & $-$
\\ \hline
18&$23$  &$0.0611$ & $1.9157 \times 10^{-8}$ & $13$ &$0.0192$
   & $4.3992 \times 10^{-8}$
\\ \hline
19&$23$  &$0.1167$ & $5.7251 \times 10^{-8}$ & $14$ &$0.0380$
   & $3.8658 \times 10^{-8}$
\\ \hline
20&$916$  &$7.1030$ & $7.8327 \times 10^{-8}$ & $5,344$ &$18.4432$
   & $3.2243 \times 10^{-8}$
\\ \hline
\end{tabular}}
}
\label{tab:18}
\end{table}

\begin{table}
\centering
\caption{
Experimental results run on the data in Table \ref{tab:50}.
``$-$'' means no convergence within $10N$
$k$-iterations.
}
\footnotesize{
{\newcommand{\q}[1]{\mc{1}{|l|}{\small\tt #1}}
\noindent
\begin{tabular}{*{8}{|c}|} \hline
\mc{2}{|c}{} &
\mc{3}{|c|}{ML($n$)BiCGStab (Alg. 5.1 in \cite{ye2012})}& \mc{3}{|c|}{ML($n$)BiCGStabt}
\\ \hline
No.
& $n$ & Iter & Time (s) & True error & Iter & Time (s) & True error
\\ \hline\hline
 1&$8$
   &$179$  &$0.2271$ & $8.4558 \times 10^{-8}$ & $184$ &$0.2283$
   & $8.1975 \times 10^{-8}$
\\ \hline
 &$16$
   &$189$  &$0.3422$ & $8.8224 \times 10^{-8}$ & $189$ &$0.3254$
   & $9.3304 \times 10^{-8}$
\\ \hline
2&$8$
   &$3,662$  &$4.0794$ & $9.6577 \times 10^{-8}$ & $2,913$ &$3.1262$
   & $9.9816 \times 10^{-8}$
\\ \hline
 &$16$
   &$2,473$  &$4.2875$ & $3.8774 \times 10^{-7}$ & $2,382$ &$3.7026$
   & $8.1546 \times 10^{-8}$
\\ \hline
3&$8$
   &$7,288$  &$8.0662$ & $6.9185 \times 10^{-8}$ & $6,503$ &$6.9870$
   & $6.0707 \times 10^{-8}$
\\ \hline
 &$16$
   &$4,108$  &$6.7744$ & $2.2397 \times 10^{-6}$ & $4,224$ &$6.5421$
   & $9.7019 \times 10^{-8}$
\\ \hline
4&$8$
   &$-$  &$-$ & $-$ & $38,054$ &$98.2040$
   & $9.6753 \times 10^{-8}$
\\ \hline
 &$16$
   &$15,895$  &$49.0668$ & $8.9919 \times 10^{-8}$ & $16,289$ &$48.1862$
   & $9.6599 \times 10^{-8}$
\\ \hline
5&$8$
   &$-$  &$-$ & $-$ & $-$ &$-$
   & $-$
\\ \hline
 &$16$
   &$24,388$  &$31.4974$ & $9.9902 \times 10^{-8}$ & $23,134$ &$27.7806$
   & $7.4456 \times 10^{-8}$
\\ \hline
6&$8$
   &$2,916$  &$5.6495$ & $8.9420 \times 10^{-8}$ & $2,900$ &$5.5693$
   & $9.0846 \times 10^{-8}$
\\ \hline
 &$16$
   &$2,670$  &$6.3544$ & $9.0912 \times 10^{-8}$ & $2,574$ &$5.9208$
   & $9.6812 \times 10^{-8}$
\\ \hline
7&$8$
   &$2,361$  &$2.1054$ & $9.5836 \times 10^{-8}$ & $2,234$ &$1.9393$
   & $5.3475 \times 10^{-8}$
\\ \hline
 &$16$
   &$2,554$  &$3.3363$ & $7.9515 \times 10^{-8}$ & $2,190$ &$2.6315$
   & $7.6172 \times 10^{-8}$
\\ \hline
8&$8$
   &$324$  &$1.4897$ & $7.6087 \times 10^{-8}$ & $327$ &$1.5177$
   & $3.1248 \times 10^{-8}$
\\ \hline
 &$16$
   &$331$  &$1.5813$ & $8.3862 \times 10^{-8}$ & $328$ &$1.5712$
   & $7.4565 \times 10^{-8}$
\\ \hline
9&$8$
   &$13,334$  &$23.6233$ & $9.7577 \times 10^{-8}$ & $9,138$ &$15.8272$
   & $9.1668 \times 10^{-8}$
\\ \hline
 &$16$
   &$1,459$  &$3.0690$ & $9.5718 \times 10^{-8}$ & $1,897$ &$3.8690$
   & $8.9255 \times 10^{-8}$
\\ \hline
10&$8$
   &$8,017$  &$17.6337$ & $1.8724 \times 10^{-5}$ & $7,934$ &$17.3321$
   & $1.5846 \times 10^{-7}$
\\ \hline
 &$16$
   &$5,417$  &$16.3895$ & $9.7412 \times 10^{-5}$ & $5,202$ &$21.2771$
   & $4.7796 \times 10^{-7}$
\\ \hline
11&$8$
   &$-$  &$-$ & $-$ & $-$ &$-$
   & $-$
\\ \hline
 &$64$
   &$-$  &$-$ & $-$ & $19,723$ &$91.8874$
   & $6.4079 \times 10^{-8}$
\\ \hline
12&$8$
   &$5,750$  &$2.4522$ & $9.2938 \times 10^{-8}$ & $6,928$ &$2.8840$
   & $8.9248 \times 10^{-8}$
\\ \hline
 &$16$
   &$4,584$  &$2.7914$ & $9.8318 \times 10^{-8}$ & $4,254$ &$2.4498$
   & $9.9897 \times 10^{-8}$
\\ \hline
13&$8$
   &$11,819$  &$43.8708$ & $9.8546 \times 10^{-8}$ & $11,807$ &$43.4537$
   & $8.1438 \times 10^{-8}$
\\ \hline
 &$16$
   &$4,440$  &$25.6544$ & $1.6295 \times 10^{-7}$ & $4,307$ &$17.8434$
   & $8.1273 \times 10^{-8}$
\\ \hline
14&$8$
   &$10$  &$0.0241$ & $4.8088 \times 10^{-8}$ & $10$ &$0.0204$
   & $4.8088 \times 10^{-8}$
\\ \hline
 &$16$
   &$10$  &$0.0221$ & $7.5853 \times 10^{-8}$ & $10$ &$0.0172$
   & $7.5853 \times 10^{-8}$
\\ \hline
15&$8$
   &$223$  &$0.2419$ & $8.4697 \times 10^{-8}$ & $227$ &$0.2368$
   & $7.6142 \times 10^{-8}$
\\ \hline
 &$16$
   &$219$  &$0.3432$ & $7.7560 \times 10^{-8}$ & $203$ &$0.3175$
   & $9.3095 \times 10^{-8}$
\\ \hline
16&$8$
   &$-$  &$-$ & $-$ & $-$ &$-$
   & $-$
\\ \hline
 &$32$
   &$1,674$  &$5.2264$ & $9.0863 \times 10^{-8}$ & $1,551$ &$5.0532$
   & $8.8965 \times 10^{-8}$
\\ \hline
17&$8$
   &$-$  &$-$ & $-$ & $-$ &$-$
   & $-$
\\ \hline
 &$64$
   &$29,202$  &$717.8035$ & $9.7427 \times 10^{-8}$ & $33,609$ &$757.6161$
   & $9.6732 \times 10^{-8}$
\\ \hline
18&$8$
   &$20$  &$0.0250$ & $5.4243 \times 10^{-8}$ & $20$ &$0.0276$
   & $5.4243 \times 10^{-8}$
\\ \hline
 &$16$
   &$21$  &$0.0296$ & $2.0159 \times 10^{-8}$ & $21$ &$0.0316$
   & $2.0159 \times 10^{-8}$
\\ \hline
19&$8$
   &$25$  &$0.0387$ & $1.6317 \times 10^{-8}$ & $25$ &$0.0511$
   & $1.6317 \times 10^{-8}$
\\ \hline
 &$16$
   &$24$  &$0.0407$ & $8.9381 \times 10^{-8}$ & $24$ &$0.0465$
   & $8.9381 \times 10^{-8}$
\\ \hline
20&$8$
   &$1,156$  &$2.8012$ & $8.8965 \times 10^{-8}$ & $1,459$ &$3.4974$
   & $8.8663 \times 10^{-8}$
\\ \hline
 &$16$
   &$617$  &$1.6847$ & $9.9117 \times 10^{-8}$ & $602$ &$1.6106$
   & $9.5743 \times 10^{-8}$
\\ \hline
\end{tabular}}
}
\label{tab:19}
\end{table}

\begin{table}
\centering
\caption{
Experimental results run on the data in Table \ref{tab:50}.
``$-$'' means no convergence within $10N$
$k$-iterations.
}
\footnotesize{
{\newcommand{\q}[1]{\mc{1}{|l|}{\small\tt #1}}
\noindent
\begin{tabular}{*{10}{|c}|} \hline
\mc{2}{|c}{} &
\mc{3}{|c|}{ML($n$)BiCGStab (Alg. 4.1 in \cite{ye2012})}& \mc{2}{|c}{} & \mc{3}{|c|}{ML($n$)BiCGStab (Alg. 4.1 in \cite{ye2012})}
\\ \hline
No.& $n$ & Iter & Time (s) & True error &No.
& $n$& Iter & Time (s) & True error
\\ \hline\hline
 1&$8$
   &$182$  &$0.2490$ & $7.0962 \times 10^{-8}$ & $11$ &$8$
   & $-$ & $-$&$-$
\\ \hline
 &$16$
   &$192$  &$0.4072$ & $8.6571 \times 10^{-8}$ & &$64$ &$-$
   & $-$ & $-$
\\ \hline
2&$8$
   &$4,110$  &$5.1773$ & $9.1479 \times 10^{-8}$ & $12$ &$8$
   & $6,495$ &$2.9844$ &$8.4379 \times 10^{-8}$
\\ \hline
 &$16$
   &$2,887$  &$5.9016$ & $8.7683 \times 10^{-8}$ & $$ &$16$
   & $5,313$ &$3.7253$ &$9.2844 \times 10^{-8}$
\\ \hline
3&$8$
   &$5,804$  &$7.2623$ & $6.2621 \times 10^{-8}$ & $13$ &$8$
   & $12,292$ &$48.8011$ &$7.4523 \times 10^{-8}$
\\ \hline
 &$16$
   &$3,464$  &$7.0040$ & $6.7311 \times 10^{-8}$ & $$ &$16$
   & $4,681$ &$30.0199$ &$7.3730 \times 10^{-8}$
\\ \hline
4&$8$
   &$34,329$  &$93.4919$ & $8.9878 \times 10^{-8}$ & $14$ &$8$
   & $13$ &$0.0291$ &$4.7338 \times 10^{-8}$
\\ \hline
 &$16$
   &$16,244$  &$54.7066$ & $8.4285 \times 10^{-8}$ & $$ &$16$
   & $14$ &$0.0308$ &$6.6326 \times 10^{-8}$
\\ \hline
5&$8$
   &$33,470$  &$30.7092$ & $7.9867 \times 10^{-8}$  & $15$ &$8$
   & $260$ &$0.3941$ &$6.1463 \times 10^{-8}$
\\ \hline
 &$16$
   &$17,147$  &$27.1564$ & $9.8394 \times 10^{-8}$ & $$ &$16$
   & $225$ &$0.5702$ &$7.4062 \times 10^{-8}$
\\ \hline
6&$8$
   &$2,915$  &$5.9639$ & $9.7400 \times 10^{-8}$ & $16$ &$8$
   & $-$ &$-$ &$-$
\\ \hline
 &$16$
   &$3,203$  &$8.4778$ & $8.2899 \times 10^{-8}$ & $$ &$32$
   & $1,330$ &$8.7912$ &$9.4115 \times 10^{-8}$
\\ \hline
7&$8$
   &$2,278$  &$2.2367$ & $5.8257 \times 10^{-8}$ & $17$ &$8$
   & $-$ & $-$&$-$
\\ \hline
 &$16$
   &$2,307$  &$3.5007$ & $7.3368 \times 10^{-8}$ & $$ &$64$
   & $28,130$ &$702.3665$ &$9.7067 \times 10^{-8}$
\\ \hline
8&$8$
   &$326$  &$1.5374$ & $4.7810 \times 10^{-8}$ & $18$ &$8$
   & $20$ &$0.0268$ &$1.7969 \times 10^{-8}$
\\ \hline
 &$16$
   &$331$  &$1.6738$ & $9.1445 \times 10^{-8}$ & $$ &$16$
   & $19$ &$0.0297$ &$8.8714 \times 10^{-8}$
\\ \hline
9&$8$
   &$12,056$  &$22.2520$ & $9.5102 \times 10^{-8}$ & $19$ &$8$
   & $24$ &$0.0399$ &$9.1399 \times 10^{-8}$
\\ \hline
 &$16$
   &$1,614$  &$3.7977$ & $9.7401 \times 10^{-8}$ & $$ &$16$
   & $24$ &$0.0458$ &$6.4031 \times 10^{-8}$
\\ \hline
10&$8$
   &$7,490$  &$17.8209$ & $1.2933 \times 10^{-7}$ & $20$ &$8$
   & $1,128$ &$2.8335$ &$9.6613 \times 10^{-8}$
\\ \hline
 &$16$
   &$5,100$  &$16.3960$ & $1.3307 \times 10^{-7}$ & $$ &$16$
   & $661$ & $1.9800$&$9.2218 \times 10^{-8}$
\\ \hline
\end{tabular}}
}
\label{tab:8-2-1}
\end{table}

\section{Concluding Remarks
}\label{con} The original motivation of developing ML($n$)BiCGStabt was to improve the stability of Algoritm 5.1 in \cite{ye2012}. From our experiments, the improvement can sometimes be significant. Since, however, the two algorithms are essentially the same in structure,
they basically share the same theoretical and numerical properties. A generalization of ML($n$)BiCGStabt to ML($n$)BiCGStabt2 and ML($n$)BiCGStabt($l$) are being carried out. They are clearly different from ML($n$)BiCGStabt in structure and thereby we expect different properties that these algorithms will have.

Now, it can be seen that ML($n$)BiCGStabt should be the first method
getting ${\bf A}^H$ involved  in its implementation in the area of
product-type or hybrid BiCG methods since the 1980 IDR method was published.



\section{Appendix} \label{sec:apen} In this section, we present a preconditioned
ML($n$)BiCGStabt algorithm together with its Matlab code.\\


\begin{algorithm}\label{alg:10-22}
{\bf ML($n$)BiCGStabt with preconditioning}
\begin{tabbing}
x\=xxx\= xxx\=xxx\=xxx\=xxx\=xxx\=xxx\=xxx\=xxx\kill
\>1. \> Choose an initial guess ${\bf
x}_0$ and $n$ vectors ${\bf q}_1, {\bf q}_2, \ldots,
{\bf q}_n$. \\
\>2. \> Compute $[{\bf f}_1, \ldots, {\bf f}_{n-1}] = {\bf M}^{-H} {\bf A}^H [{\bf q}_1, \ldots, {\bf q}_{n-1}]$, ${\bf r}_0 = {\bf b} - {\bf A} {\bf x}_0$ and ${\bf
g}_0 = {\bf r}_0$.\\
\>\> Compute $\hat{\bf g}_0 = {\bf M}^{-1}
{\bf r}_0, \,\,
 {\bf w}_0 = {\bf A}
\hat{\bf g}_0,\,\, c_0 = {\bf q}^H_{1} {\bf w}_0, \,\, e_0 = {\bf q}_1^H {\bf r}_0$. \\
\>3. \>For $j = 0, 1, 2, \ldots$ \\
\>4. \>\>For $i = 1, 2, \ldots, n-1$ \\
\>5. \>\>\>$\alpha_{jn+i} = e_{jn+i-1} / c_{jn+i-1}$;\\
\>6. \>\>\> ${\bf x}_{jn+i} = {\bf x}_{jn+i-1} + \alpha_{jn+i} \hat{\bf g}_{jn+i-1}$;
\\
\>7. \>\>\> ${\bf r}_{jn+i} = {\bf r}_{jn+i-1} - \alpha_{jn+i} {\bf w}_{jn+i-1}$;
\\
\>8. \>\>\> $e_{jn+i} = {\bf q}_{i+1}^H {\bf r}_{jn+i}$;\\
\>9. \>\>\> If $j \geq 1$\\
\>10. \>\>\>\>$\tilde{\beta}^{(jn+i)}_{(j-1)n+i} = - e_{jn+i} \big/
 c_{(j-1)n+i}$; \,\,\,\,\,\, \% $\tilde{\beta}^{(jn+i)}_{(j-1)n+i} = -\omega_{j} \beta^{(jn+i)}_{(j-1)n+i}$\\
 \>11. \>\>\>\>${\bf z}_w = {\bf r}_{jn+i} + \tilde{\beta}^{(jn+i)}_{(j-1)n+i} {\bf w}_{(j-1)n+i}$;\\
 \>12. \>\>\>\>${\bf g}_{jn+i} = \tilde{\beta}^{(jn+i)}_{(j-1)n+i} {\bf g}_{(j-1)n+i}$; \\
\>13. \>\>\>\>For $s = i+1, \ldots, n- 1$ \\
\>14. \>\>\>\>\>$\tilde{\beta}^{(jn+i)}_{(j-1)n+s} = - {\bf q}^H_{s+1}
 {\bf z}_w \big/
 c_{(j-1)n+s}$; \,\,\,\,\,\, \% $\tilde{\beta}^{(jn+i)}_{(j-1)n+s} = -\omega_{j} \beta^{(jn+i)}_{(j-1)n+s}$\\
 \>15. \>\>\>\>\>${\bf z}_w = {\bf z}_w + \tilde{\beta}^{(jn+i)}_{(j-1)n+s} {\bf w}_{(j-1)n+s}$;\\
 \>16. \>\>\>\>\>${\bf g}_{jn+i} = {\bf g}_{jn+i} + \tilde{\beta}^{(jn+i)}_{(j-1)n+s} {\bf g}_{(j-1)n+s}$; \\
\>17. \>\>\>\>End
\\
\>18. \>\>\>\>$\displaystyle{{\bf g}_{jn+i} = {\bf z}_w - \frac{1}{\omega_{j}} {\bf g}_{jn+i}}$;
\\
\>19. \>\>\>\>For $s = 0, \ldots, i - 1$ \\
\>20. \>\>\>\>\>$\beta^{(jn+i)}_{jn+s} = - {\bf f}_{s+1}^H {\bf g}_{jn+i}
\big/ c_{jn+s}$; \\
\>21. \>\>\>\>\>${\bf g}_{jn+i} = {\bf g}_{jn+i} + \beta^{(jn+i)}_{jn+s} {\bf g}_{jn+s}$; \\
\>22. \>\>\>\>End
\\
\>23. \>\>\>Else
\\
\>24. \>\>\>\>$\beta^{(jn+i)}_{jn} = - {\bf f}_{1}^H {\bf r}_{jn+i}
\big/ c_{jn}$; \\
\>25. \>\>\>\>${\bf g}_{jn+i} = {\bf r}_{jn+i} + \beta^{(jn+i)}_{jn} {\bf g}_{jn}$; \\
\>26. \>\>\>\>For $s = 1, \ldots, i - 1$ \\
\>27. \>\>\>\>\>$\beta^{(jn+i)}_{jn+s} = - {\bf f}_{s+1}^H {\bf g}_{jn+i}
\big/ c_{jn+s}$; \\
\>28. \>\>\>\>\>${\bf g}_{jn+i} = {\bf g}_{jn+i} + \beta^{(jn+i)}_{jn+s} {\bf g}_{jn+s}$; \\
\>29. \>\>\>\>End
\\
\>30. \>\>\>End
\\
\>31.\>\>\> $\hat{\bf g}_{jn+i} = {\bf M}^{-1} {\bf g}_{jn+i}
$; 
${\bf w}_{jn+i} = {\bf A} \hat{\bf g}_{jn+i}
$;\\
\>32.\>\>\> $c_{jn+i} = {\bf q}_{i+1}^H {\bf w}_{jn+i}$;\\
\>33. \>\> End\\
\>34. \>\>$\alpha_{jn+n} = e_{jn+n-1} / c_{jn+n-1}$;\\
\>35. \>\> ${\bf x}_{jn+n} = {\bf x}_{jn+n-1} + \alpha_{jn+n} \hat{\bf g}_{jn+n-1}$;
\\
\>36. \>\> $ {\bf u}_{jn+n} = {\bf r}_{jn+n-1} - \alpha_{jn+n} {\bf w}_{jn+n-1}$;
\\
\>37. \>\> $ \hat{\bf u}_{jn+n} = {\bf M}^{-1} {\bf u}_{jn+n}$;
\\
\>38. \>\> $\omega_{j+1} =  ({\bf A} \hat{\bf u}_{jn+n})^H {\bf u}_{jn+n} / \|{\bf A} \hat{\bf u}_{jn+n} \|_2^2$; \\
\>39.\>\>${\bf x}_{jn+n} = {\bf x}_{jn+n}  +\omega_{j+1} \hat{\bf u}_{jn+n}$; \\
\>40. \>\>${\bf r}_{jn+n} = -\omega_{j+1} {\bf A} \hat{\bf u}_{jn+n} +
{\bf u}_{jn+n}$; \\
\>41. \>\>$e_{jn+n} = {\bf q}_1^H {\bf r}_{jn+n}$;\\
\>42. \>\>$\tilde{\beta}^{(jn+n)}_{(j-1)n+n} = - e_{jn+n}
\big/
 c_{(j-1)n+n}$; \,\,\,\,\,\, \% $\tilde{\beta}^{(jn+n)}_{(j-1)n+n} = -\omega_{j+1} \beta^{(jn+n)}_{(j-1)n+n}$\\
 \>43. \>\>${\bf z}_w = {\bf r}_{jn+n} + \tilde{\beta}^{(jn+n)}_{(j-1)n+n} {\bf w}_{(j-1)n+n}$;\\
 \>44. \>\>${\bf g}_{jn+n} = \tilde{\beta}^{(jn+n)}_{(j-1)n+n} {\bf g}_{(j-1)n+n}$;\\
\>45. \>\>For $s = 1, \ldots, n - 1$ \\
\>46. \>\>\>$\tilde{\beta}^{(jn+n)}_{jn+s} = - {\bf q}^H_{s+1} {\bf z}_w
\big/
 c_{jn+s}$; \,\,\,\,\,\, \% $\tilde{\beta}^{(jn+n)}_{s+jn} = -\omega_{j+1} \beta^{(jn+n)}_{s+jn}$\\
 \>47. \>\>\>${\bf z}_w = {\bf z}_w + \tilde{\beta}^{(jn+n)}_{jn+s} {\bf w}_{jn+s}$;\\
 \>48. \>\>\>${\bf g}_{jn+n} = {\bf g}_{jn+n} + \tilde{\beta}^{(jn+n)}_{jn+s} {\bf g}_{jn+s}$;\\
\>49. \>\>End \\
\>50.\>\> $\displaystyle{{\bf g}_{jn+n} = {\bf z}_w  - \frac{1}{\omega_{j+1}} {\bf g}_{jn+n}}
$; 
$\hat{\bf g}_{jn+n} = {\bf M}^{-1} {\bf g}_{jn+n}
$; \\
\>51.\>\> ${\bf w}_{jn+n} = {\bf A} \hat{\bf g}_{jn+n}
$; 
$c_{jn+n} = {\bf q}_{1}^H {\bf w}_{jn+n}$;\\
\>52. \> End
\end{tabbing}
\end{algorithm}

\vspace{.2cm}

{\bf
Matlab code of Algorithm \ref{alg:10-22}}
\begin{tabbing}
x\=xxx\=
xxx\=xxx\=xxx\=xxx\=xxx\=xxx\=xxx\=xxx\=xxx\=xxx\=xxx\=xxx\kill \>1.
\>function $[x,err,iter,flag] = mlbicgstabt(A,x,b,Q,M,max\_it,tol, kappa)$\\
\>2.\\
\>3.\>\% input:\>\>\>$A$:\,\,\, N-by-N matrix. $M$: \,\,\,N-by-N preconditioner matrix. \\
\>4.\>\% \>\>\>$Q$: \, N-by-n shadow matrix $[{\bf
q}_1,\ldots,{\bf q}_n]$. $x$: initial guess.\\
\>5.\>\% \>\>\> $b$:\,\,\, right hand side vector. $max\_it$:\,\,\, maximum number
of iterations.\\
\>6.\>\%\>\>\>$tol$:\,\,\, error tolerance.\\
\>7.\>\%\>\>\>$kappa$:\>\>\> (real number) minimization step controller:\\
\>8.\>\%\>\>\>\>\>\> $kappa = 0$, standard minimization\\
\>9.\>\%\>\>\>\>\>\> $kappa > 0$, Sleijpen-van der Vorst minimization \\
\>10.\>\% output:\>\>\>$x$: solution computed. $err$: error norm. $iter$: number of iterations
performed.\\
\>11.\>\% \>\>\>$flag$:\>\>\> $= 0$, solution found to
tolerance\\
\>12.\>\% \>\>\>\>\>\>  $= 1$, no convergence given $max\_it$ iterations\\
\>13.\>\% \>\>\>\>\>\>$= -1$, breakdown. \\
\>14.\>\% storage:\>\>\> $F$: $N \times (n-1)$ matrix. $G, Q, W$: $N \times n$ matrices. $A, M$: $N \times N$ matrices.\\
\>15.\>\% \>\>\>$x, r, g\_h, z, b$: $N \times 1$ matrices. $c$: $1 \times n$ matrix.\\
\>16. \\
\>17.\>\>$N = size(A,2); \,\,n = size(Q,2)$;\\
\>18.\>\>$G = zeros(N,n);\,\, W = zeros(N,n)$;\,\,\,\,\, \%
initialize
work spaces\\
\>19.\>\>if $n > 1$, $F = zeros(N,n-1)$; end\\
\>20.\>\>$c = zeros(1,n)$;\,\,\,\,\,\,\,\,\,\,\,\,\, \% end initialization\\
\>21.\>\>\\
\>22.\>\>$iter = 0;\,\, flag = 1;\,\, bnrm2 = norm(b)$;\\
\>23.\>\>if $bnrm2 == 0.0$,\, $bnrm2 = 1.0$;\, end\\
\>24.\>\>$r = b - A*x; \,\, err = norm( r ) / bnrm2$;\\
\>25.\>\>if $err < tol$,\, $flag = 0$;\,\, return,\, end\\
\>26.\>\> \\
\>27.\>\> if $n > 1$, $F = M' \backslash (A' * Q(:,1: n-1))$; end\\
\>28.\>\>$G(:,1) = r;\,\, g\_h = M \backslash r; \,\, W(:,1) = A*g\_h; \,\, c(1) = Q(:,1)'*W(:,1)$;\\
\>29.\>\>if $c(1) == 0$,\, $flag = -1$;\,\, return,\, end \\
\>30.\>\>$e = Q(:,1)'*r$; \\
\>31.\>\>\\
\>32.\>\>  for $j = 0:max\_it$\\
\>33.\>\>\>for $i = 1:n-1$\\
\>34.\>\>\>\> $alpha = e / c(i)$;
\, $x = x + alpha*g\_h$; \,
$r = r - alpha*W(:, i)$;\\
\>35.\>\>\>\> $err = norm(r)/bnrm2$;\,
$iter = iter + 1$;\\
\>36.\>\>\>\>if $err < tol$,\, $flag = 0$;\,\, return,\, end\\
\>37.\>\>\>\>if $iter >= max\_it$,\, return,\, end\\
\>38.\>\>\>\>\\
\>39.\>\>\>\>$e = Q(:,i+1)'*r$; \\
\>40.\>\>\>\>if $j >= 1$\\
\>41.\>\>\>\>\>$beta = -e / c(i+1)$; \\
\>42.\>\>\>\>\>$W(:,i+1) = r + beta*W(:,i+1)$; \\
\>43.\>\>\>\>\>$G(:,i+1) = beta*G(:,i+1)$; \\
\>44.\>\>\>\>\> for $s = i+1 : n-1$\\
\>45.\>\>\>\>\>\>$beta = -Q(:,s+1)'*W(:,i+1)/c(s+1)$; \\
\>46.\>\>\>\>\>\>$W(:,i+1) = W(:,i+1) + beta*W(:,s+1)$;\\
\>47.\>\>\>\>\>\>$G(:,i+1) = G(:,i+1) + beta*G(:,s+1)$;\\
\>48.\>\>\>\>\>end \\
\>49.\>\>\>\>\>$G(:,i+1) = W(:,i+1) - G(:,i+1)./omega$; \\
\>50.\>\>\>\>\>for $s = 0:i-1$\\
\>51.\>\>\>\>\>\>$beta = -F(:,s+1)'*G(:,i+1) / c(s+1)$;\\
\>52.\>\>\>\>\>\>$G(:,i+1) = G(:,i+1) + beta*G(:,s+1)$; \\
\>53.\>\>\>\>\>end \\
\>54.\>\>\>\>else \\
\>55.\>\>\>\>\>$beta = -F(:,1)'*r / c(1)$; \,
$G(:,i+1) = r + beta*G(:,1)$;\\
\>56.\>\>\>\>\>for $s = 1:i-1$\\
\>57.\>\>\>\>\>\>$beta = -F(:,s+1)'*G(:,i+1) / c(s+1)$;\\
\>58.\>\>\>\>\>\>$G(:,i+1) = G(:,i+1) + beta*G(:,s+1)$; \\
\>59.\>\>\>\>\>end \\
\>60.\>\>\>\>end \\
\>61.\>\>\>\>$g\_h = M \backslash G(:,i+1)$; \,
$W(:,i+1) = A*g\_h$; \\
\>62.\>\>\>\>$c(i+1) = Q(:,i+1)'*W(:,i+1)$;\\
\>63.\>\>\>\>if $c(i+1) == 0$,\, $flag = -1$;\,\, return,\, end \\
\>64.\>\>\>end \\
\>65.\>\>\> $alpha = e / c(n)$;\,
$x = x + alpha*g\_h$; \,
$r = r - alpha*W(:,n)$; \\
\>66.\>\>\>
$err = norm(r)/bnrm2$;\\
\>67.\>\>\>if $err < tol$,\, $flag = 0; \,\, iter = iter + 1$;\, return,\, end \\
\>68.\>\>\>$g\_h = M \backslash r; \,\, z = A*g\_h$;\, $omega = z'*z$; \\
\>69.\>\>\>if $omega == 0$,\, $flag = -1$;\, return,\, end\\
\>70.\>\>\>$rho = z'*r$;\,
$omega =  rho / omega$; \\
\>71.\>\>\>if $kappa > 0$\\
\>72.\>\>\>\>$rho = rho / (norm(z)*norm(r))$;\,
$abs\_om = abs(rho)$;\\
\>73.\>\>\>\>if ($abs\_om < kappa$) \& ($abs\_om \sim = 0$)\\
\>74.\>\>\>\>\>$omega = omega*kappa/abs\_om$;\\
\>75.\>\>\>\>end\\
\>76.\>\>\>end\\
\>77.\>\>\>if $omega == 0$,\, $flag = -1$;\, return,\, end\\
\>78.\>\>\>$x = x + omega*g\_h$;\,
$r = r - omega*z$;\\
\>79.\>\>\>$err = norm(r)/bnrm2$;\,
$iter = iter + 1$;\\
\>80.\>\>\>if $err < tol$,\, $ flag = 0$;\,\, return,\, end\\
\>81.\>\>\>if $iter >= max\_it$,\, return,\, end \\
\>82.\>\>\> \\
\>83.\>\>\>$e = Q(:,1)'*r; \,\, beta = - e/c(1)$; \\
\>84.\>\>\>$W(:,1) = r + beta*W(:,1)$;\,
$G(:,1) = beta*G(:,1)$; \\
\>85.\>\>\> for $s = 1:n-1$ \\
\>86.\>\>\>\>$beta = -Q(:,s+1)'*W(:,1)/c(s+1)$;\\
\>87.\>\>\>\>$W(:,1) = W(:,1) + beta*W(:,s+1)$; \\
\>88.\>\>\>\>$G(:,1) = G(:,1) + beta*G(:,s+1)$;\\
\>89.\>\>\> end \\
\>90.\>\>\>$G(:,1) = W(:,1) - G(:,1)./omega$;\,
$g\_h = M \backslash G(:,1)$; \\
\>91.\>\>\>$W(:,1) = A * g\_h$; \,
$c(1) = Q(:,1)'*W(:,1)$;\\
\>92.\>\>\>if $c(1) == 0$,\, $flag = -1$;\,\, return,\, end \\
\>93.\>\>  end
\end{tabbing}

\vspace{.2cm}


\begin{thebibliography}{10}







\bibitem{dsyyz} {\sc L. Du, T. Sogabe, B. Yu, Y. Yamamoto, S.-L. Zhang},
{\it A block IDR($s$) method for nonsymmetric linear systems with multiple right-hand sides},
J. Comput. Appl. Math. 235(2011), no. 14, 4095-4106.













\bibitem{fletcher}
{\sc R. Fletcher}, {\it Conjugate gradient methods for indefinite
systems,} volume 506 of Lecture Notes Math., pages 73-89.
Springer-Verlag, Berlin-Heidelberg-New York, 1976.



\bibitem{fgn}
{\sc R. Freund, M. Gutknecht and N. Nachtigal}, {\it An
implementation of the look-ahead Lanczos algorithm for non-Hermitian
matrices}, SIAM J. Sci. Comput., 14(1993), pp. 137-158.





\bibitem{gs2013}
{\sc M. Gijzen and P. Sonneveld}, {\it Algorithm 913: an elegant IDR($s$) variant that efficiently exploits bi-orthogonality properties}, ACM Trans. Math. Software, 38(2011), pp. 5:1--5:19.

\bibitem{gut30}
{\sc M. H. Gutknecht}, {\it A completed theory of the unsymmetric
Lanczos process and related algorithms. Part I.}, SIAM J. Matrix
Anal. Appl., 13(1992), pp.594-639.


\bibitem{gut}
------, {\it Variants of BICGStab for matrices with complex spectrum,}
SIAM J. Sci. Comput., 14 (1993), pp.\ 1020--1033.


\bibitem{gut31}
------, {\it A completed theory of the unsymmetric Lanczos process and related
algorithms. Part II.}, SIAM J. Matrix Anal. Appl. 1994, 15:15-58.

\bibitem{gut10}
------, {\it Lanczos-type solvers for nonsymmetric linear systems of equations,}
Acta Numerica, 6 (1997), pp. 271-397.

\bibitem{gut1}
------, {\it IDR Explained,}
ETNA 36 (2010), 126--148.

\bibitem{gutzem}
{\sc Martin H. Gutknecht and Jens-Peter M. Zemke},
{\it Eigenvalue computations based on IDR}, Bericht 145, TUHH, Institute of Numerical Simulation, May 2010.






\bibitem{joubert1}
{\sc W. D. Joubert}, {\it Generalized conjugate gradient and Lanczos
methods for the solution of nonsymmetric systems of linear
equations}, Ph.D. thesis and Tech. Report CNA-238, Center for
Numerical Analysis, University of Texas, Austin, TX, 1990.

\bibitem{joubert2}
------, {\it Lanczos methods for the solution of nonsymmetric systems of
linear equations}, SIAM Journal on Matrix Analysis and Applications
1992; 13:926-943.

\bibitem{loher}
{\sc D. Loher}, {\it Reliable nonsymmetric block Lanczos algorithms}, Ph.D. thesis,
Swiss Federal Institute of Technology, Zurich, 2006.

\bibitem{neum}
{\sc A. Neumaier}, {\it Oral presentation at the Oberwolfach meeting ``Numerical Linear Algebra''}, Oberwolfach, April 1994.
















\bibitem{parlett}
{\sc B. N. Parlett, D. R. Taylor, and Z. A. Liu}, {\it A look-ahead
Lanczos algorithm for unsymmetric matrices}, Math. Comp., 44(1985), pp.105-124.




\bibitem{saad82}
{\sc Y. Saad}, {\it The Lanczos biorthogonalization algorithm and
other oblique projection methods for solving large unsymmetric
systems}, SIAM Journal on Numerical Analysis, 19(1982), pp.485-506.


\bibitem{saad2}
{\sc Y. Saad and M. H. Schultz}, {\it GMRES: A generalized minimal
residual algorithm for solving nonsymmetric linear systems,} SIAM J.
Sci. Statist. Comput., 7 (1986), pp.\ 856--869.






\bibitem{SF}
{\sc G. L. G. Sleijpen and D. R. Fokkema}, {\it BiCGSTAB($l$) for
linear equations involving unsymmetric matrices with complex
spectrum,} ETNA, 1:11-32, 1993.


\bibitem{sv2013}
{\sc G.L.G. Sleijpen and M. B. van Gijzen}, {\it
Exploiting BiCGstab($l$) strategies to induce dimension reduction}, SIAM J. Sci. Comput. 32(2010), no. 5, 2687-2709.










\bibitem{Svan}
{\sc G. L. G. Sleijpen and H. A. van der Vorst},
{\it Reliable
updated residuals in hybrid Bi-CG methods,} Computing, 56 (1996),
pp. 141--163.




\bibitem{sonn}
{\sc P. Sonneveld}, {\em CGS, a fast Lanczos-type solver for
nonsymmetric linear systems,} SIAM J. Sci. Statist. Comput., 10
(1989), pp.\ 36--52.


\bibitem{sonn2013}
{\sc P. Sonneveld}, {\em On the convergence behavior of IDR($s$) and related methods},
SIAM J. Sci. Comput., 34(5), A2576-A2598.


\bibitem{sonn1}
{\sc P. Sonneveld and M. van Gijzen}, {\em IDR(s): a family of
simple and fast algorithms for solving large nonsymmetric systems of linear
equations,} SIAM J. Sci. Comput. 31(2008), no. 2, pp. 1035-1062.

\bibitem{tasu}
{\sc M. Tanio and M. Sugihara}, {\it GBi-CGSTAB($s,l$): IDR($s$) with higher-order stabilization polynomials}, J. Comput. Appl. Math. 235(2010), no. 3, 765-784.

\bibitem{van}
{\sc H. A. van der Vorst}, {\it Bi-CGSTAB: A fast and smoothly
converging variant of Bi-CG for the solution of nonsymmetric linear
systems,} SIAM J. Sci. Statist. Comput., 12 (1992), pp.\ 631--644.





































\bibitem{ws}
{\sc P. Wesseling and P. Sonneveld}, {\em Numerical experiments with
a multiple grid and a preconditioned Lanczos type method}, Lecture Notes in Mathematics,
vol. 771, pp.543-562,
Springer Verlag, Berlin, Heidelberg, New York, 1980.

\bibitem{ye2012}
{\sc M. Yeung}, {\em ML($n$)BiCGStab: Refomulation, Analysis and Implementation},
Numer. Math. Theor. Meth. Appl. 5 (2012), pp. 447-492.


\bibitem{yeung7}
{\sc M. Yeung}, {\em An introduction to ML($n$)BiCGStab}, available
at http://arxiv.org/abs/1106.3678.
Proceedings of Boundary Elements and Other Mesh Reduction Methods XXXIV, edited by Brebbia \& Popov, 2012, WITpress.

\bibitem{yeungboley}
{\sc M. Yeung and D. Boley}, {\em Transpose-free multiple Lanczos
and its application in Pad\'e approximation}, Journal of
Computational and Applied Mathematics, Vol 177/1 pp. 101-127, 2005.


\bibitem{yeungchan}
{\sc M. Yeung and T. Chan}, {\em ML($k$)BiCGSTAB: A BiCGSTAB variant
based on multiple Lanczos starting vectors}, SIAM J. Sci. Comput.,
Vol. 21, No. 4, pp.~1263-1290, 1999.



\bibitem{zhang}
{\sc Shao-Liang Zhang}, {\em GPBi-CG: Generalized product-type
methods based on Bi-CG for solving nonsymmetric linear systems},
SIAM J. Sci. Comput., 18:537-551, 1997.

\end{thebibliography}
\end{document}